\documentclass[11pt,a4paper,reqno]{article}
\usepackage{amsfonts,color, amsmath,amssymb}
\usepackage{tikz}

\input{xy}
\xyoption{all}

\renewcommand{\mod}{\textup{mod}\,}

\newcommand{\supp}{\textup{supp}\,}
\newcommand{\Hom}{\mbox{\rm Hom}}

\newcommand{\Ext}{\textup{Ext}}

\newcommand{\End}{\textup{End}}

\newcommand{\add}{\textup{add}\,}

\newcommand{\C}{\mathcal{C}}

\newcommand{\id}{\mbox{$\mathbb I$}}

\newtheorem{thm}{Theorem}[section]
\newtheorem{prop}[thm]{Proposition}
\newtheorem{cor}[thm]{Corollary}
\newtheorem{lem}[thm]{Lemma}
\newtheorem{defn}{Definition}
\newtheorem{rem}[thm]{Remark}

\addtocounter{section}{0}
\begin{document}


\title{Projective dimension of modules over cluster-tilted algebras}

 \author{Louis Beaudet \thanks{Partially supported by the NSERC of Canada and the University of Sherbrooke} 
 \and Thomas Br\"ustle \thanks{Partially supported by the NSERC of Canada, University of Sherbrooke and Bishop's University} 
  \and   Gordana Todorov \thanks  {Partially supported by the NSA-grant MSPF-08G-228}}





\maketitle

\begin{abstract}
We study the projective dimension of finitely generated modules over cluster-tilted algebras $\End_{\mathcal C}(T)$ where $T$ is a cluster-tilting object in a cluster category $\C$.
It is well-known that all $\End_{\mathcal C}(T)$-modules are of the form $\Hom_{\mathcal C}(T,M)$ for some object $M$ in $\C$, and since $\End_{\mathcal C}(T)$  is Gorenstein of dimension 1, the projective dimension of
$\Hom_{\mathcal C}(T,M)$ is either zero, one or infinity.
We define in this article the ideal $I_M$ of $\End_\C (T[1])$ given by all endomorphisms that factor through $M$, and show that the $\End_{\mathcal C}(T)$-module $\Hom_{\mathcal C}(T,M)$ has infinite projective dimension precisely  when $I_M$ is non-zero.
Moreover, we apply the results above to characterize the location of modules of infinite projective dimension in the Auslander-Reiten quiver of cluster-tilted algebras of type $A$ and $D$.
\end{abstract}


\begin{section}{Introduction}\label{sect intro}
Cluster categories for acyclic quivers $Q$ were introduced in  \cite{BMRRT}, and, for type
$A_n$ also in \cite{CCS}, as a means for a better understanding of the
cluster algebras of Fomin and Zelevinsky \cite{FZ1,FZ2}. They are
defined as orbit categories  ${\mathcal C}_Q = D^b( \mod kQ) / \tau^{-1}[1]$ of the  bounded derived
categories of finitely generated $kQ$-modules.

In \cite{BMR1},  Buan, Marsh and Reiten defined the cluster-tilted
algebras as follows. Let 
$T$ be a cluster-tilting object in $\mathcal{C}_Q$, that is,  an object
such that 
$\Ext^1_{{\mathcal C}_Q}(T,T)=0$ and 
the number of non-isomorphic indecomposable summands of
$T$ equals the number of vertices of the quiver $Q$.
Then the endomorphism algebra  $\End_{{\mathcal C}_Q}(T)$ is called a
{\it  cluster-tilted algebra}. 
Since then, these algebras have been the subject of many
investigations, see, for instance, \cite{ABS1,ABS2,ABS3,BMR1,BMR3,BRS,CCS,CCS2,KR}.

While cluster categories serve as our main motivation, our main result holds in a more general situation, namely for any triangulated category $\C$ with a maximal 1-orthogonal object $T$ (see Definition \ref{def}).
Generalizing results obtained for cluster-tilted algebras and 2-Calabi-Yau categories in \cite{BMR1,KR}, 
Koenig and Zhu show in \cite{KZ} that the functor $\Hom_{\C}(T,-)$ induces an equivalence $\C/\add T[1] \to \mod (\End_{\C} T)$.
Moreover,  the algebra $\End_{\C}(T)$ is Gorenstein of dimension 1, thus we know all $\End_{\C}(T)$-modules have projective dimension zero, one or infinity.
In this paper we characterize the modules of infinite projective dimension.

By the above equivalence, indecomposable $\End_{\C}(T)$-modules correspond to indecomposable objects $M$ in $\C$ which do not belong to $\add T[1]$. For any such object, we denote by $I_M$ the ideal  of $\End_\C (T[1])$ given by all endomorphisms that factor through $M$ and call it {\it factorization ideal of M}.

Our main theorem is the following:

\begin{thm} Let $\mathcal C$ be a triangulated category with a maximal 1-orthogonal object $T$. Let $M$ be an indecomposable object in $\mathcal C$ which does not belong to  ${\rm add } \, T[1]$.
 Then the $\End_{\C}(T)$-module $\Hom_{\C}(T,M)$ has infinite projective dimension precisely when the factorization ideal $I_M$ is non-zero.
\end{thm}

We discuss in the last section of the paper the fact  that the indecomposable $\End_{\C}(T)$-modules of infinite projective dimension are lying on certain hammock-like subquivers of the Auslander-Reiten quiver of $\C$ provided $\C=\C_{Q}$ where $Q$ is Dynkin, and we give a precise description for quivers of type $A_n$ and $D_n$.

We would like to thank Kiyoshi Igusa and Alex Lasnier for helpful discussions.
\end{section}


\begin{section}{Preliminaries}\label{sect 1}
\begin{defn}\label{def} Let $\C$ be a triangulated category.
An object $T$ of $\C$ is called  maximal 1-orthogonal if  \add T is contravariantly finite and covariantly finite, and an object 
$ X \in \C$ belongs to $ \add T$ if and only if $\Ext_{\C}^1(X,T)=0=\Ext_{\C}^1(T,X).$
\end{defn}
Denote by $C$ the endomorphism ring $C:= \End_{\C} T$, and by  $\mod C$  the category of finitely presented $C-$modules. We recall from \cite{KZ} that the functor $\Hom_{\C}(T,-)$ induces an equivalence $\C/\add T[1] \to \mod C$, and
moreover, under this equivalence the projective $C-$modules are those of the form $\Hom_{\C}(T,T^0)$ where $T^0$ belongs to $\add T$.

For an object $M$ in $\C$, define $I_M$ to be the set of all $f \in \End_{\C} \, T[1]$ such that there exist morphisms $ g:T[1] \to M$  and $h:M\to T[1] $  such that $f = h \circ g $.
Clearly, $I_M$ is an ideal of the algebra $\End_{\C} \, T[1] \cong C$ and it is easy to see that $I_{M\oplus N} = I_M + I_N$.
Since our focus is on the quotient $\C/\add T[1] \simeq \mod C$, we will usually consider only objects $M$ in $\C$ such  that $\add M \cap \add T[1] =0$.

We recall some results that we need later on:

\begin{lem}[\cite{KR}]\label{KR}
Let $T$ be a maximal 1-orthogonal object in a triangulated category $\C$, and let
\[
\xymatrix@C=30pt{Y[-1] \ar[r]& X \ar[r] & T^0
  \ar[r]^f & Y}   \]
  be a triangle in $\C$ with $f: T^0 \to Y$  a contravariant $\add T$-approximation of $Y$. Then $X $ belongs to $\add T$.
\end{lem}

\begin{lem}\label{KZ} Let $\C$ be a triangulated category with a maximal 1-orthogonal  object $T$. Let $g : X \to Y$ be a morphism in $\C$ which is a part of a triangle
\[
\xymatrix@C=30pt{Z[-1] \ar[r]^{ h }& X \ar[r]^{g} & Y
  \ar[r]^f & Z}.   \]
 Then $\Hom_{\C}(T,g)$ is a monomorphism in $\mod C$ if and only if $\Hom_{\C}(T,h)=0$.
 \end{lem}
 
 {\bf Proof.}
This follows directly by applying the functor $\Hom_{\C}(T,-)$ to the given triangle .
   \hfill $\Box$ \bigskip
 
When there is no danger of confusion, we will sometimes instead of
\[
\xymatrix@C=50pt{\Hom_\C(T,M) \ar[r]^{ \Hom_\C(T,f) }&\Hom_\C(T,N)}
 \]
 write one of the following simplified forms:
 \[
\xymatrix@C=40pt{\Hom(T,M) \ar[r]^{ (T,f) }&\Hom(T,N)}
 \]
\[
\xymatrix@C=40pt{(T,M) \ar[r]^{ (T,f) }&(T,N)}.
 \]

 The following theorem has been shown in \cite{KR} in the context of 2-Calabi-Yau categories, and in the general context of triangulated categories in \cite{KZ}:
 \begin{thm}[{\cite[Cor 4.5]{KZ}}]\label{Gorenstein}
Let $\C$ be a triangulated category and $T$ a maximal 1-orthogonal object in $\C$. Then the endomorphism ring $C$ of $T$ is Gorenstein of Gorenstein dimension at most one.
 \end{thm}

 \begin{rem}\label{pdim}
 As it has been observed in \cite{KR}, being Gorenstein of Gorenstein dimension at most one implies that each $C$-module has either projective dimension at most one, or is of infinite projective dimension.
 \end{rem}

\end{section}


\begin{section}{Proof of the main theorem}\label{sect 3}
Throughout this section, we fix an indecomposable object $M$ in $\C$ which is not a direct summand of $T[1]$.
The aim is to prove that Hom$_{\mathcal C}(T,M)$ has infinite projective dimension if and only if the factorization ideal $I_M\neq 0$. We first show the following implication:
\begin{subsection}{The factorization ideal $I_M$ of $M$ is non-zero when $C$-module $\Hom_{\C}(T,M)$ has infinite projective dimension.}
{\bf Proof.}
Assume $\Hom_{\C}(T,M)$ has infinite projective dimension, and consider a projective cover
\[
\xymatrix@C=30pt{\Hom_{\C}(T,T^0) \ar[r]^{ (T,f) }&\Hom_{\C}(T,M)
  \ar[r] & 0}   \]

in $\mod C$ defined by a morphism $f:T^0 \to M$ in $\C$.
Choose a triangle in $\C$
\[
\xymatrix@C=30pt{M[-1] \ar[r]^{ \beta[-1] }& T^1 \ar[r] & T^0
  \ar[r]^f & M}   \]
containing $f$, then by Lemma \ref{KR} we know that $T^1$ belongs to $\add T$.
Applying the functor $(T,\ ):=\Hom_{\C}(T,-)$ to that triangle yields an exact sequence in $\mod C$
\[
\xymatrix@C=40pt{(T,M[-1]) \ar[r]^{ (T,\beta[-1]) }& (T,T^1) \ar[r] & (T,T^0)   \ar[r] & (T,M)   \ar[r] & (T,T^1[1])} \]
where $\Hom_{\C}(T,T^1[1])=0$ since $T$ is a maximal 1-orthogonal object.
This implies that the morphism $\Hom_{\C}(T,\beta[-1])$ is non-zero, since otherwise the projective dimension of $\Hom_{\C}(T,M)$ would be at most one.
Choose a morphism $\alpha[-1]$ in $\Hom_{\C}(T,M[-1])$ whose image under $\Hom_{\C}(T,\beta[-1])$ is non-zero, that is, the composition
\[
\xymatrix@C=30pt{T \ar[r]^{ \alpha[-1] }&M[-1] \ar[r]^{\beta[-1]} & T^1}
 \]
 is non-zero.
 This yields the non-zero composition
 \[
\xymatrix@C=30pt{T[1] \ar[r]^{ \alpha }&M \ar[r]^{\beta} & T^1[1]}
 \]
As $T^1[1]$ is a non-trivial summand in $\add T[1]$, we conclude that there is a non-zero element in the factorization ideal $I_M$ of $M$.
\hfill $\Box$ \bigskip
\end{subsection}

\begin{subsection}{The factorization ideal $I_M$ of $M$ is non-zero only if $\Hom_{\C}(T,M)$ has infinite projective dimension.}
{\bf Proof.}
By remark \ref{pdim}, a $C$-module $\Hom_{\C}(T,M)$ of finite projective dimension has projective dimension zero or one. It suffices to show that in both cases the ideal $I_M$ is zero.

{\bf Case 1:} Assume $\Hom_{\C}(T,M)$ has projective dimension $0$.\medskip\\
Then $M$ belongs to $\add T$, and every composition of morphisms
 \[
\xymatrix@C=30pt{T[1] \ar[r]^{ \alpha }&M \ar[r]^{\beta} & T[1]}
 \]
must be zero since $\beta \in \Hom_{\C}(M,T[1])=0$.
Therefore $I_M=0$ in this case.

{\bf Case 2:} Assume $\Hom_{\C}(T,M)$ has projective dimension $1$.\medskip\\
Thus there is a projective resolution
\[
\xymatrix@C=30pt{(1) & 0 \ar[r]^{  }& (T,T^1) \ar[r]^{(T,g)} & (T,T^0)   \ar[r] & (T,M)   \ar[r] & 0}
\]
in $\mod C$.
We choose a triangle
\[
\xymatrix@C=30pt{(2)\quad &    & Z[-1] \ar[r]^{ h }& T^1 \ar[r]^{g} & T^0
  \ar[r]^f & Z &}   \]
  in $\C$, and by lemma \ref{KZ} we conclude that $\Hom_{\C}(T,h)=0$. Since moreover $\Hom_{\C}(T,T^1[1])=0$, applying $(T,\ )=\Hom_{\C}(T,-)$ to the triangle (2) yields a short exact sequence in $\mod C$

\[
\xymatrix@C=30pt{(3) & 0 \ar[r] & (T,T^1) \ar[r]^{(T,g)} & (T,T^0) \ar[r]^{(T,f)} & (T,Z)  \ar[r] & 0}.
\]
Since the two short exact sequences (1) and (3) start with the same morphism $(T,g)$, we conclude that their cokernels $\Hom_{\C}(T,M)$ and $\Hom_{\C}(T,Z)$ are isomorphic.
This implies in the category $\C$ that the objects $M$ and $Z$ differ only by summands in $ \add T[1]$. But we assumed $M$ to be indecomposable and not isomorphic to an object in $\add T[1]$, hence $M$ is isomorphic to a summand of $Z$. We denote by $\iota: M \to Z$ the corresponding section with retraction $\rho:Z \to M$.

Given a factorization
 \[
\xymatrix@C=30pt{T[1] \ar[r]^{ \alpha }&M \ar[r]^{\beta} & T[1]}
 \]

the aim is to show that $\beta \alpha =0$ and hence $I_M=0$.
We consider the composed maps
 \[
\xymatrix@C=30pt{T[1] \ar[r]^{ \iota \alpha }&Z \ar[r]^{\beta\rho} & T[1]}
 \]

 and insert them into a commutative diagram formed using the triangle (2):
 \[
\xymatrix@C=30pt{ & & T \ar[r]\ar[d]^x & 0 \ar[d]\ar[r] & T[1] \ar[d]^{\iota\alpha} \ar[r]^{\id} & T[1] \ar[d]^{x[1]} \\
(4) \qquad & Z[-1] \ar[r]^{ h }& T^1 \ar[d]^y\ar[r]^{g} & T^0  \ar[d]\ar[r]^f & Z \ar[d]^{\beta\rho}\ar[r] & T^1[1]\ar[d]^{y[1]}\\
& & T \ar[r] & 0 \ar[r] & T[1]  \ar[r]^{\id} & T[1] }
  \]
The existence of a morphism $x$ making the upper part of the diagram commutative is guaranteed by the following commutative square between two triangles in $\C$:
\[
\xymatrix@C=30pt{  0 \ar[d]\ar[r] & T[1] \ar[d]^{\iota\alpha}  \\
 T^0  \ar[r]^f & Z}
  \]

Likewise, the existence of a morphism $y$ making the lower part of the diagram commutative is guaranteed by the following square
\[
\xymatrix@C=30pt{  T^0 \ar[d]\ar[r]^f & Z \ar[d]^{\beta\rho}  \\
 0  \ar[r] & T[1]}
  \]

 which commutes since $\Hom_{\C}(T^0,T[1])=0$.
 From the commutative square
 \[
\xymatrix@C=30pt{  T \ar[d]^x \ar[r] & 0 \ar[d]  \\
 T^1  \ar[r]^g & T^0}
  \]

 in (4) we conclude that $gx=0$, and hence $x=0$ since $\Hom_{\C}(T,g)$ is a monomorphism.
 This implies $x[1]=0$ and therefore, by the commutativity of (4),
 \[
 0 = y[1] x[1] = \beta\rho\iota\alpha = \beta\alpha
 \]
 which implies $I_M=0$.

  \hfill $\Box$ \bigskip

  \begin{rem}
\rm  We can apply the dual proof to obtain that the $\End \,_{\C} T-$module $\Hom_{\C}(T,M)$ has infinite injective dimension precisely when the factorization ideal $I_M$ is non-zero.
However, it is easy to see for any Gorenstein algebra that the modules of infinite projective dimension are exactly the modules of infinite injective dimension, so it is no surprise they satisfy the same condition.
  \end{rem}
\begin{rem}
 \rm We showed in the proof of case 2 in 3.2 that there exists the triangle (2) in $\C$ where  $M$ and $Z$ differ only by summands in $ \add T[1]$, that is $Z\cong M \oplus Z'$ where $Z' \in \add T[1]$. We then continue to show in (4) that there is no morphism from $T[1]$ to $T[1]$ factoring through $Z$. 
 This implies that $Z'=0$, hence we can extract from the proof above the following lifting property:
  \end{rem}

\begin{cor}
Every short exact sequence
\[
\xymatrix@C=30pt{(1) & 0 \ar[r]^{  }& (T,T^1) \ar[r]^{(T,g)} & (T,T^0)   \ar[r]^{(T,f)} & (T,M)   \ar[r] & 0}
\]
in $\mod C$ with $T^0,T^1 \in \add T$ can be lifted to 
 a triangle
\[
\xymatrix@C=30pt{(2)\quad &    & M[-1] \ar[r]^{ h }& T^1 \ar[r]^{g} & T^0
  \ar[r]^f & M &}   \]
  in $\C$.
\end{cor}
  \end{subsection}

\end{section}


\begin{section}{Hammocks, swings and rays}\label{ham}

We discuss in this section how to compute the objects $M \in \C$ with $I_M \neq 0$ and give a detailed description for the modules of infinite projective dimension for cluster algebras of Dynkin type $A$ and $D$.
We assume throughout that $\C$ is a triangulated category with  a maximal 1-orthogonal object $T=T_1 \oplus \cdots \oplus T_n$.
 Of course one has $I_M \neq 0$ if and only if there are indices $i,j$ such that there exists a non-zero morphism between the indecomposable summands $T_i[1] $ and $ T_j[1] $ of $T[1]$ factoring through $M$.
Define the following full subquiver
$$H(i,j) = \{ X \in \C | \mbox{ there is } 0 \neq \lambda:T_i[1] \to T_j[1] \mbox{ factoring through } X \}$$
of the Auslander-Reiten quiver $\Gamma(\C)$ of $\C$.
Our main theorem implies that the functor $\Hom_{\C}(T,-)$ induces a bijection between the set 
$$ \bigcup_{i,j} (H(i,j) \backslash T[1] )$$
 and the indecomposable $C-$modules of infinite projective dimension.
 Thus, to determine all indecomposable $C-$modules of infinite projective dimension, it is sufficient to compute all the full subquivers $H(i,j)$ of the Auslander-Reiten quiver $\Gamma(\C)$. 
 We illustrate this procedure with the following example.
\bigskip

{\bf Example:} Let $C$ be the cluster-tilted algebra given by the quiver

\begin{displaymath}
\xymatrix @C=15pt @R=20pt @!R @!C {  & 5\ar[rd]^b & & &3\ar[ld]_d & \\
          &  & 4 \ar[r]^c & 1\ar[rr]_e & & 2\ar[lu]_f \\
          &  6\ar[ru]^a & &  & & \\}
\end{displaymath}

with relations $de=ef=fd=0$. 
The algebra $C$ is cluster-tilted of type ${\mathbb D}_6$, given as endomorphism ring of a cluster-tilting object $T=T_1 \oplus \cdots \oplus T_6$ of the cluster category $\C$ of type ${\mathbb D}_6$.
The full translation subquivers $H(i,j)$ of the Auslander-Reiten quiver $\Gamma(\C)$ can be calculated easily using mesh-relations and starting from one indecomposable $T_i[1]$.

It turns out that only for the three pairs $(i,j) = (2,1),(1,3),(3,2)$ the set $H(i,j) \backslash T[1]$ is non-empty (there are no $C-$modules lying on a path from, say, $T_4[1]$ to $T_6[1]$).
The modules in $H(3,1)$ are indicated in Figure 1 in blue, whereas the modules in  $H(2,1)$ and $H(1,3)$ are red.

{\tiny
\begin{displaymath}
\xymatrix @C=4pt @R=0.1pt @!R @!C { &&&&&&&&&&&&&&&\\          &&                   &                  &                           &  T_{6}[1] \ar[rdd]      &                  &  P_{6} \ar[rdd]        &  & \color{blue}\txt{5\\4\\1} \ar[rdd]&  & \txt{6\\43\\1} \ar[rdd]& & & & &\\&
                     &                   &                  &                           &  T_{5}[1] \ar[rd]       &                  & P_{5} \ar[rd]  & & \color{blue}\txt{6\\4\\1} \ar[rd]& & \txt{5\\43\\1} \ar[rd]& &  &  &&\\&
                     &                   &      .~.~.            & T_{4}[1] \ar[rd]\ar[ruu]\ar[ru]  &               & P_{4}\ar[rd]\ar[ru]\ar[ruu]        &  & \color{blue}\txt{65\\$4^{2}$\\$1^{2}$\\2} \ar[rd]\ar[ru]\ar[ruu] &  & \color{blue}\txt{65\\$4^{2}3$\\$1^{2}$} \ar[rd]\ar[ru]\ar[ruu]& & \txt{65\\$4^{2}3$\\1} \ar[rd] & .~.~.& &&\\&
                     &                   & T_{1}[1] \ar[rd]\ar[ru] &                           &  P_{1} \ar[rd]\ar[ru]&                  & \color{blue}\txt{4\\1}\ar[rd]\ar[ru] &  & \txt{65\\$4^{2}3$\\$1^{2}$\\2} \ar[rd]\ar[ru] & & \color{blue}\txt{65\\$4^{2}$\\1} \ar[rd]\ar[ru]& &  \txt{65\\$4 3$\\1} \ar[rd]&  &\\&
                     &  \color{red}\txt{3} \ar[dr]\ar[ur] &                  & \color{red}\txt{2} \ar[rd]\ar[ru]          &                  & \color{blue}\txt{1} \ar[rd]\ar[ru] &  & \txt{43\\1} \ar[rd]\ar[ru] & & \txt{65\\$4^{2}$\\1\\2} \ar[rd]\ar[ru]& &\color{blue}\txt{65\\4\\1} \ar[rd]\ar[ru] & & \txt{3} \ar[rd]&\\ \ar@{--}@*{[red]}'[ururur]'[rrrrrr]'[drrrrr]'[rrru]'[dr]'[]&
          T_{2}[1] \ar[ru]  &       & P_{2} \ar[ru]        &  \ar@{--}@*{[blue]}'[urururur]'[urururururu]'[ururururrr]'[rrrrrrrrrr]'[drrrrrrrrr]'[rrrrruuu]'[dr]'[]                         & T_{3}[1] \ar[ru]        &                  & P_{3}\ar[ru] & & \txt{4} \ar[ru]& & \txt{65\\4\\1\\2} \ar[ru]& & T_{2}[1] \ar[ru] & & P_{2}&\\  & &&&&&&&&&&&&&&&\\}
            \end{displaymath}
}
\begin{center} Figure 1
\end{center}

\begin{subsection}{Hammocks}

For each $1 \le i \le n$ we define the {\em left hammock} $H_i$ to be the full subquiver of the Auslander-Reiten quiver $\Gamma(\C)$ given by
\[ H_i = \supp \Hom_\C(T_i[1],-)
\]

The {\em right hammock} $_iH$ is defined as the full subquiver of the Auslander-Reiten quiver $\Gamma(\C)$ given by
\[ _iH = \supp \Hom_\C(-,T_i[1])
\]

These left and right hammocks can be computed easily using the mesh relations in $\Gamma(\C)$. Hammocks in cluster categories are described in \cite{Ri} where the right hammock $H_i$ is given as the set of vertices where the cluster-hammock function $h_{T_i[2]}$ takes positive values, and in case $\C$ is of Dynkin type, the left hammock $_iH$ coincides with the set of vertices where the cluster-hammock function $h_{T_i}$ takes positive values.

In case $H(i,j) \neq \emptyset$ it follows from the definition that $T_j[1]$ lies in the left hammock $H_i$ and that $T_i[1]$ lies in the right hammock $_jH$.
Moreover, it is clear that $H(i,j)$ is a subset
\[ H(i,j) \subseteq (H_i \cap\,  _jH)
\]
However, it is not clear when this inclusion is proper. In order to investigate $H(i,j)$ in more detail, we first study the possible relative positions of $T_i[1]$ and $T_j[1]$ in the Auslander-Reiten quiver $\Gamma(\C)$ provided $\C$ is of Dynkin type.
Since $T$ is maximal 1-orthogonal, the possible location of $T_j[1]$ in the left hammock $H_i$ is quite restricted: As $T_j[1]$ is not contained in the shifted hammock $H_i[-1] = \supp \Ext^1_\C(T_i[1],-)$, it must be on the left rim $H_i \backslash H_i[-1]$ of $H_i$ (these vertices are referred to as the projective vertices of the translation quiver $H_i$).

Each vertex lying on a sectional path starting in $T_i[1]$ is a projective vertex of the translation quiver $H_i$, we first investigate this situation:

\begin{lem}\label{sectional}
Let $\C$ be of Dynkin type. Assume that there is a sectional path $\sigma: T_i[1] \to \cdots \to T_j[1]$ in the Auslander-Reiten quiver $\Gamma(\C)$.
Then $H(i,j)$ is formed by all vertices on $\sigma$ and
 \[ H(i,j) = (H_i \cap\,  _jH)
\]
 \end{lem}
{\bf Proof.}
Clearly every vertex on the sectional path $\sigma$ belongs to $H(i,j)$. 
Conversely, let $X \in (H_i \cap\,  _jH)$, thus there is a path 
$$\gamma: T_i[1] \to \cdots \to X \to \cdots \to T_j[1]$$ in  $\Gamma(\C)$.
Since $\sigma$ is sectional, it is not possible that $\sigma$ and $\gamma$ differ by mesh relations. And since $\C$ is of Dynkin type, there is only one sectional path from $T_i[1] $ to $T_j[1]$. Thus $\sigma = \gamma$ and therefore $H(i,j) = (H_i \cap\,  _jH)$.
  \hfill $\Box$ \bigskip
  
\begin{subsection}{Type $A$}
Assume that $\C$ is a cluster category of type  $A$. Then the left hammock $H_i$ is of rectangular shape with unique source $T_i[1]$ as depicted in Figure 2. The support of $\Ext^1_\C(T_i[1],-)$ is given by the dark shaded rectangle, and therefore $T_j[1]$ must lie on a sectional path
\[ T_i[1] \to \cdots \to T_j[1]
\]
provided $H(i,j)$ is non-empty.
The description of $ H(i,j)$  is given in Lemma \ref{sectional}.

\end{subsection}
$$
\begin{tikzpicture}
[inner sep=0.5mm]
\filldraw[fill=gray!20!white, draw=black] (0,0) -- (3,3) -- (5,1) -- (2,-2) -- cycle;
   \filldraw[fill=gray!50!white, draw=black]
    (1,0) -- (4,3) -- (6,1) -- (3,-2) -- cycle;
\draw   (1,0) -- (4,3) -- (6,1) -- (3,-2) -- cycle;
\draw (3,3)  -- (3.5,2.5);
 \draw [dashed] (-3,3) -- (8,3);
 \draw [dashed] (3.5,2.5) -- (5,1) -- (2,-2);
\draw [dashed] (-3,-2) -- (8,-2);
 \draw [ultra thick,blue] (0,0) -- (2,2);
 \draw [ultra thick] (2,2) -- (3,3);
 \draw [ultra thick] (0,0) -- (2,-2);
 \node at (-0.5,0) {$T_i[1]$};
 \node at (0,0) [shape=circle,draw,fill] {};
 \node at (1.5,2) {$T_j[1]$};
 \node at (2,2) [shape=circle,draw,fill] {};
 \end{tikzpicture}
$$
\begin{center} Figure 2
\end{center}
\end{subsection}

\begin{subsection}{Type $D$}
We suppose now that $\C$ is a cluster category of type $D$. The Auslander-Reiten quiver of $\C$ differs in this case from type $A$ by a chain of meshes with three middle terms, as indicated at the bottom of Figure 3. We first consider the case when $T_i[1]$ is not located at the top or bottom $\tau-$orbits, that is we assume that $\Hom(T_i[1],T_i) \neq 0$.
 As shown in Figure 3, in this case the projective vertices of the left hammock $H_i$ lie in two connected components, where we indicated in each copy one possible location for $T_j[1]$. The first component is given by all points living on a sectional path starting in $T_i[1]$, and  if $T_j[1]$ is one of these points, then $H(i,j)$ is described by Lemma \ref{sectional}.

$$
\begin{tikzpicture}
[inner sep=0.5mm]
 
 \filldraw[fill=gray!50!white, draw=black] (0,0) -- (2,2) -- (4.5,-0.5) -- (7,2) -- (8.5,0.5) -- (6,-2) -- (5.5,-2.5) -- (5,-2) -- (4.5,-2.5) -- (4,-2) -- (3.5,-2.5) -- (3,-2) -- (2.5,-2.5) -- cycle;
\filldraw[fill=gray!20!white, draw=black] (-1,0) -- (1,2) -- (3.5,-0.5) -- (6,2) -- (7.5,0.5) -- (5,-2) -- (4.5,-2.5) -- (4,-2) -- (3.5,-2.5) -- (3,-2) -- (2.5,-2.5) -- (2,-2) -- (1.5,-2.5) -- cycle;
  \draw[dashed] (0,0) -- (2,2) -- (4.5,-0.5) -- (7,2) -- (8.5,0.5) -- (6,-2) -- (5.5,-2.5) -- (5,-2) -- (4.5,-2.5) -- (4,-2) -- (3.5,-2.5) -- (3,-2) -- (2.5,-2.5) -- cycle;
\draw [dashed] (-3,2) -- (9,2);
 \draw (-3,-2) -- (9,-2);
 \node at (-1.5,0) {$T_i[1]$};

 \foreach \i in {-3,-2.5,...,9}
    \node at (\i,-2) [inner sep=0.3mm,shape=circle,draw,fill] {};
 \foreach \i in {-3,-2,...,8}
  \draw[xshift=\i cm] (0,-2) -- (0.5,-2.5);
\foreach \i in {-3,-2,...,8}
  \draw[xshift=\i cm] (1,-2) -- (0.5,-2.5);
 \foreach \i in {-3,-2,...,8}
  \draw[xshift=\i cm] (0,-2) -- (0.5,-1.5);
\foreach \i in {-3,-2,...,8}
  \draw[xshift=\i cm] (1,-2) -- (0.5,-1.5);
\draw [ultra thick] (1.5,-2.5) -- (-1,0) -- (1,2);
\draw [ultra thick] (1.5,-2) -- (1,-2);
\draw [ultra thick] (6,2) -- (4,0);
 \node at (-1,0) [shape=circle,draw,fill] {};
\node at (5,1) [shape=circle,draw,fill] {};
\node at (0,1) [shape=circle,draw,fill] {};
\node at (4.3,1) {$T_j[1]$};
\node at (-0.7,1) {$T_j[1]$};

 \end{tikzpicture}
$$
\begin{center} Figure 3
\end{center}
We therefore assume now that $T_j[1]$ does not lie on a sectional path starting in $T_i[1]$, and the aim is to describe $H(i,j)$ in this case. 
It is clear from Figure 3 that there are precisely two vertices $X,Y$ of  the Auslander-Reiten quiver $\Gamma(\C)$ which are middle terms of a mesh of width 3 and such that there is a sectional path from $T_i[1]$ to $X$ and $Y$ and a sectional path from $X$ and $Y$ to $T_j[1]$. We call the collection of all vertices on such a sectional path a {\em swing with fixpoints $T_i[1],T_j[1]$} and denote it by $S(T_i[1],T_j[1])$. An example is depicted in Figure 4.

$$
\begin{tikzpicture}
[inner sep=0.5mm]
 
\draw [dashed] (-3,2) -- (9,2);
 \draw (-3,-2) -- (9,-2);
 
 \foreach \i in {-3,-2.5,...,9}
    \node at (\i,-2) [inner sep=0.3mm,shape=circle,draw,fill] {};
 \foreach \i in {-3,-2,...,8}
  \draw[xshift=\i cm] (0,-2) -- (0.5,-2.5);
\foreach \i in {-3,-2,...,8}
  \draw[xshift=\i cm] (1,-2) -- (0.5,-2.5);
 \foreach \i in {-3,-2,...,8}
  \draw[xshift=\i cm] (0,-2) -- (0.5,-1.5);
\foreach \i in {-3,-2,...,8}
  \draw[xshift=\i cm] (1,-2) -- (0.5,-1.5);
\draw [ultra thick] (1.5,-2.5) -- (-1,0);
\draw [ultra thick] (2,-2) -- (1,-2);
\draw [ultra thick] (5,1) -- (1.5,-2.5);
 \node at (-1,0) [shape=circle,draw,fill] {};
\node at (5,1) [shape=circle,draw,fill] {};
\node at (-1.5,0) {$T_i[1]$};
\node at (5.7,1) {$T_j[1]$};
\node at (1.5,-2) [shape=circle,draw,fill] {};
\node at (1.5,-2.5) [shape=circle,draw,fill] {};
\node at (1.5,-1.7) {$X$};
\node at (1.5,-2.8) {$Y$};

 \end{tikzpicture}
$$
\begin{center} Figure 4
\end{center}

\begin{prop}\label{swing}
Let $\C$ be of Dynkin type $D$. Assume that $T_i[1]$ and $T_j[1]$ are located in the Auslander-Reiten quiver $\Gamma(\C)$ in such a way that $H(i,j) \neq 0$,  $\Hom(T_i[1],T_i) \neq 0$ and such that there is no sectional path from $T_i[1]$ to $T_j[1]$.
Then $H(i,j)$ is a swing with fixpoints $T_i[1],T_j[1]$ as described in Figure 4. 
Moreover,  \[ H(i,j) \neq (H_i \cap\,  _jH)
\]

in case $S(T_i[1],T_j[1]) \neq (H_i \cap\,  _jH)$.
 \end{prop}
{\bf Proof.}
It is clear from the mesh relations in $\Gamma(\C)$ that all vertices on the swing $S(T_i[1],T_j[1])$ belong to $H(i,j)$.
Since $H(i,j) \subseteq (H_i \cap\,  _jH)$ we only have to show that every vertex in $ (H_i \cap\,  _jH)\,  \backslash \, S(T_i[1],T_j[1])$ does not belong to $H(i,j)$.
We draw in Figure 5 the set $H_i \cup \,_jH$ in light gray, and its intersection $H_i \cap \,_jH$ is indicated in dark gray.
We denote by $Z_1, Z_2, \ldots$ the vertices lying on the sectional path from $X$ or $Y$ to $T_j[1]$. 
Note that $\dim \Hom(T_i[1],Z_1)=2$ but $\dim \Hom(T_i[1],T_j[1])=1$, so the dimension drops to 1 along the sectional path. We denote by $Z_t$ the vertex such that $\dim \Hom(T_i[1],Z_t)=2$ but $\dim \Hom(T_i[1],Z_{t+1})=1$.

$$
\begin{tikzpicture}
[inner sep=0.5mm]
 
 \filldraw[fill=gray!20!white, draw=white] (0,-2) -- (0.5,-2.5) -- (1,-2) -- (1.5,-2.5) -- (5,1) -- (4,2) -- (0.5,-1.5) -- (-3,2) -- (-3.5,1.5) -- cycle;
 \filldraw[fill=gray!20!white, draw=white] (-1,0) -- (1,2) -- (3.5,-0.5) -- (6,2) -- (7.5,0.5) -- (5,-2) -- (4.5,-2.5) -- (4,-2) -- (3.5,-2.5) -- (3,-2) -- (2.5,-2.5) -- (2,-2) -- (1.5,-2.5) -- cycle;
\draw [dashed] (-3,2) -- (9,2);
 \draw (-3,-2) -- (9,-2);
 \draw (5,-2) -- (7.5,0.5) -- (6,2) -- (5,1) -- (4,2) -- (2.5,0.5) -- (1,2) -- (-1,0) -- (-3,2) -- (-3.5,1.5) -- (0.5,-2.5) ;
 \filldraw[fill=gray!50!white, draw=gray!20!white] (1.5,-2.5) -- (3.5,-0.5) -- (2.5,0.5) -- (0.5,-1.5) -- cycle;
 \draw [dashed] (1.5,-2.5) -- (3.5,-0.5) -- (2.5,0.5) -- (0.5,-1.5) -- cycle;
 \foreach \i in {-3,-2.5,...,9}
    \node at (\i,-2) [inner sep=0.3mm,shape=circle,draw,fill] {};
 \foreach \i in {-3,-2,...,8}
  \draw[xshift=\i cm] (0,-2) -- (0.5,-2.5);
\foreach \i in {-3,-2,...,8}
  \draw[xshift=\i cm] (1,-2) -- (0.5,-2.5);
 \foreach \i in {-3,-2,...,8}
  \draw[xshift=\i cm] (0,-2) -- (0.5,-1.5);
\foreach \i in {-3,-2,...,8}
  \draw[xshift=\i cm] (1,-2) -- (0.5,-1.5);
\draw [ultra thick] (1.5,-2.5) -- (-1,0);
\draw [ultra thick] (2,-2) -- (1,-2);
\draw [ultra thick] (5,1) -- (1.5,-2.5);
 \node at (-1,0) [shape=circle,draw,fill] {};
\node at (5,1) [shape=circle,draw,fill] {};
\node at (-1.5,0) {$T_i[1]$};
\node at (5.7,1) {$T_j[1]$};
\node at (1.5,-2) [shape=circle,draw,fill] {};
\node at (1.5,-2.5) [shape=circle,draw,fill] {};
\node at (1.5,-1.7) {$X$};
\node at (1.5,-2.8) {$Y$};
\node at (2,-2.4) {$Z_1$};
\node at (3.7,-0.8) {$Z_t$};
\node at (3.4,-1.3) {$Z_{t-1}$};
\node at (3.0,0.3) {$W_t$};
\node at (3.8,0.8) {$W_{t+1}$};
\node at (4.4,-0.3) {$Z_{t+1}$};
\node at (2,-2) [shape=circle,draw,fill] {};
\node at (2.5,-1.5) [shape=circle,draw,fill] {};
\node at (3,-1) [shape=circle,draw,fill] {};
\node at (3,0) [shape=circle,draw,fill] {};
\node at (3.5,0.5) [shape=circle,draw,fill] {};
\node at (3.5,-0.5) [shape=circle,draw,fill] {};
\node at (4,0) [shape=circle,draw,fill] {};
\node at (4.5,0.5) [shape=circle,draw,fill] {};

 \end{tikzpicture}
$$
\begin{center} Figure 5
\end{center}

The vertex $Z_t$ has two immediate predecessors in $\Gamma(\C)$, denoted $Z_{t-1}$ and $W_t$ in Figure 5.
Assume that $W_t$ belongs to $(H_i \cap\,  _jH)$, thus $\dim \Hom(T_i[1],W_t)=1=\dim \Hom(W_t,T_j[1])$ and we denote by $\alpha: T_i[1] \to W_t$ and $\beta: W_t \to T_j[1]$ non-zero morphisms.
By the mesh relations, each morphism from $W_t$ to $T_j[1]$ factors through $W_{t+1}$, thus we can write $\beta$ as a composition $\beta = \gamma \circ \delta$ where $\delta$ is a morphism from $W_t$ to  $W_{t+1}$ and $\gamma$ is a morphism from  $W_{t+1}$ to $T_j[1]$.
But by the definition of $Z_t$ we know that $W_{t+1} \not\in H_i$ and thus $\delta \circ \alpha = 0$.
this shows that $\beta \circ \alpha = \gamma \circ \delta \circ \alpha = 0$ and therefore we conclude that $W_t$ does not belong to $H(i,j)$.
An analogous argument shows that none of the predecessors of $W_t$ except those in the swing $S(T_i[1],T_j[1])$ belongs to $H(i,j)$.
  \hfill $\Box$ \bigskip

 The last case to be considered is when $T_i[1]$ is located at the top or bottom $\tau-$orbits, that is when $\Hom(T_i[1],T_i) = 0$, and
 by symmetry we can assume the same for $T_j[1]$.

\begin{prop}
Let $\C$ be of Dynkin type $D$. Assume that $T_i[1]$ and $T_j[1]$ are located in  the Auslander-Reiten quiver $\Gamma(\C)$ in 
such a way that $H(i,j) \neq 0$,  $\Hom(T_i[1],T_i) \neq 0$ and $\Hom(T_j[1],T_j) \neq 0$.
Then  \[ H(i,j) = H_i \cap \,  _jH.
\]
 \end{prop}

{\bf Proof.}
 In case $T_i[1]$ and $T_j[1]$ are both located at the top $\tau-$orbit it is easy to see that 
  $H(i,j)$  is a swing with fixpoints $T_i[1],T_j[1]$, and in this case one also has $H(i,j) = H_i \, = \,  _jH$.
  In case one of $T_i[1]$ and $T_j[1]$ is located at the top $\tau-$orbit and the other at the two bottom orbits, then there must be a sectional path from $T_i[1]$ to $T_j[1]$ and the situation is covered in Lemma \ref{sectional}.
  Therefore we focus now on the situation when both $T_i[1]$ and  $T_j[1]$ are located at the two bottom $\tau-$orbits.
The situation is described schematically in Figure 6: The alternating white spots on the bottom two $\tau-$orbits indicate vertices $X$ where $\Hom(T_i[1],X)=0$. Then $T_j[1]$ must be, as a projective vertex in $H_i$, located on a position $X[1]$. The vertices $X$ appear one in each mesh of three middle terms, which turns the mesh relations into a commutativity. Thus the hammocks are easily computed as in case $\mathbb A$, and one clearly sees that $H(i,j) = H_i \cap \,  _jH$, which is given in Figure 6 by the shaded region.
 \hfill $\Box$ \bigskip

$$
\begin{tikzpicture}
[inner sep=0.5mm]
 
 \filldraw[fill=gray!20!white, draw=black] (2.5,0.5)  -- (5,-2) -- (4.5,-2.5) -- (4,-2) -- (3.5,-2.5) -- (3,-2) -- (2.5,-2.5) -- (2,-2) -- (1.5,-2.5) -- (1,-2)  -- (0.5,-2.5) --(0,-2) -- cycle;
\draw [dashed] (-3,2) -- (9,2);
 \draw (-3,-2) -- (9,-2);
  \foreach \i in {-3,-2.5,...,9}
    \node at (\i,-2) [inner sep=0.3mm,shape=circle,draw,fill] {};
 \foreach \i in {-3,-2,...,8}
  \draw[xshift=\i cm] (0,-2) -- (0.5,-2.5);
\foreach \i in {-3,-2,...,8}
  \draw[xshift=\i cm] (1,-2) -- (0.5,-2.5);
 \foreach \i in {-3,-2,...,8}
  \draw[xshift=\i cm] (0,-2) -- (0.5,-1.5);
\foreach \i in {-3,-2,...,8}
  \draw[xshift=\i cm] (1,-2) -- (0.5,-1.5);
 \node at (-0.5,-2.5) [shape=circle,draw,fill] {};
\node at (5.5,-2.5) [shape=circle,draw,fill] {};
\node at (0.5,-2.5) [inner sep=0.9mm,shape=circle,draw=gray,fill=white] {};
\node at (4.5,-2.5) [inner sep=0.9mm,shape=circle,draw=gray,fill=white] {};
\node at (2.5,-2.5) [inner sep=0.9mm,shape=circle,draw=gray,fill=white] {};
\node at (1.5,-2.0) [inner sep=0.9mm,shape=circle,draw=gray,fill=white] {};
\node at (3.5,-2.0) [inner sep=0.9mm,shape=circle,draw=gray,fill=white] {};
\node at (-0.5,-3) {$T_i[1]$};
\node at (5.7,-3.0) {$T_j[1]$};

 \end{tikzpicture}
$$
\begin{center} Figure 6
\end{center}

\end{subsection}
\end{section}

\end{document}